\documentclass{article}
\title{Linear relations of zeroes of the zeta-function}
\author{D. G. Best\footnote{Supported by NSERC CGS-M and the University of Lethbridge SGS Fellowship}\\Department of Mathematics and Computer Science\\University of Lethbridge, AB T1K 3M4, Canada \\ darcy.best@uleth.ca \and T. S. Trudgian\footnote{Supported by ARC Grant DE120100173}\\Mathematical Sciences Institute\\ The Australian National University, ACT 0200, Australia\\ timothy.trudgian@anu.edu.au}

\usepackage{url}
\usepackage{amsthm}
\usepackage{amsmath}
\usepackage{amssymb}
\usepackage{graphicx,float}
\usepackage{algorithm}
\usepackage{booktabs}
\newtheorem{Lem}{Lemma}

\newtheorem{theorem}{Theorem}
\newtheorem{defn}{Definition}

\numberwithin{equation}{section}

\DeclareGraphicsExtensions{.jpg,.png}

\begin{document}
\maketitle

\begin{abstract}
This article considers linear relations between the non-trivial zeroes of the Riemann zeta-function. The main application is an alternative disproof to Mertens' conjecture by showing that $\limsup_{x\rightarrow\infty} M(x) x^{-1/2} \geq 1.6383$, and  $\liminf_{x\rightarrow\infty} M(x) x^{-1/2} \leq -1.6383.$
\end{abstract}

\section{Introduction and Results}
It is not known whether any non-trivial zeroes of the zeta-function are linearly dependent over the rationals. That is, no one has found an $N\geq 1$ and integers $c_{1}, \ldots, c_{N}$, not all zero, for which
\begin{equation}\label{lindep}
\sum_{n=1}^{N} c_{n} \gamma_{n} = 0,
\end{equation}
where $\rho_{n} = \beta_{n} + i\gamma_{n}$ is the $n$th non-trivial zero of $\zeta(s)$. It seems that Ingham \cite{Ingham1942} was the first to consider (\ref{lindep}). His paper concerned, \textit{inter alia}, Mertens' conjecture that $M(x) = \sum_{n\leq x} \mu(n) \leq x^{\frac{1}{2}}$, where $\mu(n)$ is the M\"{o}bius function. Ingham showed that Mertens' conjecture implies that there are infinitely many linear dependencies as given in (\ref{lindep}). Since there seems to be no intrinsic reason why (\ref{lindep}) should be true, Ingham expressed doubts about Mertens' conjecture. Indeed, Mertens' conjecture was shown to be false by Odlyzko and te Riele in \cite{MertDis}.

Ingham's result is `doubly infinite': there are infinitely many choices for the $c_{n}$ and infinitely many $N$. Bateman et al.\ \cite{Bateman} proved a `singly infinite' result: Mertens' conjecture implies that there are infinitely many sums of the type (\ref{lindep}) in which the $c_{n}$ are integers, not all zero, $|c_{n}|\leq 2$, and at most one $c_{n} = \pm 2.$ Furthermore, they considered all of these permissible sums for $1\leq N \leq 20$ and showed that no linear dependencies exist. We extend their table in Appendix A.

The contrapositive to the statement given by Bateman et al.\ is an interesting one: if there is no relation of the type (\ref{lindep}) for $|c_{n}| \leq 2$, then Mertens' conjecture is false. This singly infinite result was reduced to a finite result with the work of Grosswald \cite{Grosswald}. Following Grosswald, we are able to prove the following

\begin{theorem}\label{Ours}
\begin{equation*}
\limsup_{x\rightarrow\infty} M(x) x^{-1/2} \geq 1.6383, \quad \liminf_{x\rightarrow\infty} M(x) x^{-1/2} \leq -1.6383.
\end{equation*}
\end{theorem}

This improves on the result of Kotnik and te Riele \cite{Kotte} who showed that $\limsup_{x\rightarrow\infty} M(x) x^{-1/2} \geq 1.218$ and $\liminf_{x\rightarrow\infty} M(x) x^{-1/2} \leq -1.229$. An added feature to the approach in this paper is that the bounds given in Theorem \ref{Ours} for the $\limsup$ and the $\liminf$ are equal. 

\section{Outline}

Suppose $g(x)$ is a piecewise-continuous real function that is bounded on finite intervals. Suppose also that
\begin{equation}\label{bigg}
G(s) = \int_{1}^{\infty} g(x) x^{-s-1}\, dx
\end{equation}
originally absolutely convergent for $\sigma > \sigma_{a}$, say, can be continued analytically to $\sigma = \sigma_{0}$. Moreover, assume that one can write the principal part of $G(s)$ as
\begin{equation}\label{bigg2}
H(s) = G(s) - \left( \frac{r_{0}}{s- \sigma_{0}} + \sum_{\gamma} \frac{r_{\gamma}}{s- (\sigma_{0} + i\gamma)}\right),
\end{equation}
where $H(s)$ is analytic for $s= \sigma_{0} + it$, where $|t|<T$. Here, $\gamma$ is an element of some finite set of numbers $0 < |\gamma| <T$. Ingham \cite[Thm 1]{Ingham1942} proved\footnote{Ingham actually proved a slightly different version from that given above --- see \cite{HaselP} and \cite[p.\ 86]{AndStark} for details.} that for any $T>0$ and any $x_{0}$,
\begin{equation}\label{InghamThm1}
\liminf_{x\rightarrow\infty} \frac{g(x)}{x^{\sigma_{0}}} \leq r_{0} + \sum_{-T < \gamma < T} r_{\gamma} \left(1- \frac{|\gamma|}{T}\right) x_{0}^{i\gamma} \leq \limsup_{x\rightarrow\infty} \frac{g(x)}{x^{\sigma_{0}}}.
\end{equation}

To exhibit large negative values of $\liminf_{x\rightarrow\infty} \frac{g(x)}{x^{\sigma_{0}}}$, for example, one hopes to align the arguments of the complex terms $x_{0}^{i\gamma}$ so that they all pull in the same direction. If the numbers $\gamma$ are independent one can achieve this using Kronecker's theorem.

Grosswald's idea is to extract some partial information by weakening the hypothesis that the $\gamma$s are linearly independent. This weaker version of linear independence has also been considered by Anderson and Stark \cite{AndStark} and Diamond \cite{Diamond}. The following version is that given in \cite{AndStark}.

Let $T>0$, and let $\Gamma$ denote a set of positive numbers. Define $\Gamma'$ as a subset of $\Gamma$ such that every $\gamma\in\Gamma'$ lies in the range $0< \gamma< T$. Finally, let $\{N_{\gamma}\}$ be a set of positive integers defined for $\gamma \in \Gamma'$.

\begin{defn}\label{defcon1}
The elements of $\Gamma'$ are $\{N_{\gamma}\}$-independent in $\Gamma \cap [0, T]$ if
\begin{equation}\label{ind1}
\sum_{\gamma\in\Gamma'} c_{\gamma} \gamma = 0, \quad \textrm{with } |c_{\gamma}| \leq N_{\gamma},
\end{equation}
implies that all $c_{\gamma} = 0$ and for any $\gamma^{*}\in \Gamma \cap [0, T]$,
\begin{equation}\label{ind2}
\sum_{\gamma\in\Gamma'} c_{\gamma} \gamma = \gamma^{*}, \quad \textrm{with } |c_{\gamma}| \leq N_{\gamma},
\end{equation}
implies that $\gamma^{*} \in \Gamma'$, that $c_{\gamma*} = 1$, and that all other $c_{\gamma} = 0$.
\end{defn}

With this definition, it is possible to prove the following diluted version of (\ref{InghamThm1}).
\begin{theorem}[Anderson and Stark]\label{AndThm}
If the elements of $\Gamma'$ are $\{N_{\gamma}\}$-independent in $\Gamma \cap[0,T]$, then
\begin{equation}\label{AndS1}
\liminf_{x\rightarrow\infty} \frac{g(x)}{x^{\sigma_{0}}} \leq r_{0} - \sum_{\gamma \in \Gamma'} \frac{2N_{\gamma}}{N_{\gamma} + 1} |r_{\gamma}| \left(1- \frac{\gamma}{T}\right)
\end{equation}
and
\begin{equation}\label{AndS2}
\limsup_{x\rightarrow\infty} \frac{g(x)}{x^{\sigma_{0}}} \geq r_{0} + \sum_{\gamma \in \Gamma'} \frac{2N_{\gamma}}{N_{\gamma} + 1} |r_{\gamma}| \left(1- \frac{\gamma}{T}\right).
\end{equation}
\end{theorem}

\subsection{Mertens' Conjecture}
Ingham \cite[pp.\ 318-319]{Ingham1942}  gave a proof of the classical result that
\begin{equation}\label{LastIn}
\liminf_{x\rightarrow \infty} \frac{M(x)}{\sqrt{x}} = -\infty, \quad\quad \limsup_{x\rightarrow \infty} \frac{M(x)}{\sqrt{x}} = \infty.
\end{equation}
follows if either the Riemann hypothesis is false or not all the zeroes are simple.
Henceforth we assume the Riemann hypothesis and the simplicity of the zeroes.
In (\ref{bigg}), take $g(x) = M(x)$ and so $\sigma_{0} = \frac{1}{2}$ and, in (\ref{bigg2}), $r_{0} = 0$ and $r_{\gamma} = (\rho\zeta'(\rho))^{-1}.$ Here $\rho = \frac{1}{2} + i\gamma$ is a typical non-trivial zero of the zeta-function.

Ingham \cite[Thm.\ A]{Ingham1942} used (\ref{InghamThm1}) to show that if the zeroes $\gamma$ are linearly independent, then (\ref{LastIn}) is true. We now follow Grosswald's approach.

 Let $M \geq 1$ be given. Choose $\Gamma' = \{\gamma_{1}, \ldots, \gamma_{M}\}$ and choose $T = \gamma_{L+1}-\epsilon$ for some small $\epsilon$, where $L \geq M$. Thus, $\Gamma\cap [0,T] = \{\gamma_{1}, \ldots, \gamma_{L}\}$. Therefore, provided that (\ref{ind1}) and (\ref{ind2}) are satisfied, we have by Theorem \ref{AndThm}
\begin{equation}\label{AndS3}
\liminf_{x\rightarrow\infty} \frac{M(x)}{x^{1/2}} \leq - \sum_{n=1}^{M} \frac{2N_{\gamma_{n}}}{N_{\gamma_{n}} + 1} \frac{1}{|\rho_{n} \zeta'(\rho_{n})|} \left(1- \frac{\gamma_{n}}{T}\right)
\end{equation}
and
\begin{equation}\label{AndS4}
\limsup_{x\rightarrow\infty} \frac{M(x)}{x^{1/2}} \geq  \sum_{n=1}^{M} \frac{2N_{\gamma_{n}}}{N_{\gamma_{n}} + 1} \frac{1}{|\rho_{n} \zeta'(\rho_{n})|} \left(1- \frac{\gamma_{n}}{T}\right).
\end{equation}

\section{Computation}\label{comp}

Previous disproofs of Mertens' conjecture have utilized the basis reduction algorithm first described by  Lenstra, Lenstra and Lov\'asz in \cite{lll}, called LLL-reduction. We also employ the use of this robust algorithm, but in a different way. In order to explain our process, we shall first provide the algorithm so that it may be used as a road map while reading this section. Although we will be directing our attention at zeroes, all of these processes work with any set of real numbers. We will also assume that all $N_{\gamma_{n}}$s are equal, and we denote their value by $N_{\gamma}$.

\begin{algorithm}[H]\label{main-algo}
\textbf{Input:} $(k,n,T) \in \mathbb{Z} \times \mathbb{Z} \times \mathbb{R}$.\\
\textbf{Output:} $N_\gamma$.
\begin{enumerate}
 \item Define $K := 10^k$ and $m := \left|\Gamma \cap [0,T]\right|$.
 \item Compute all zeroes, $\gamma$, such that $0 <\gamma \leq T$ to at least $k$ decimal digits of precision.
 \item Sort the elements of $\Gamma \cap [0,T]$ using $\prec$ in (\ref{heavy-order}). Label the {\it heaviest} element $\gamma_1$, the next {\it heaviest} element $\gamma_2$, $\ldots$ , and the least {\it heavy} element $\gamma_m$.
 \item Define $\Gamma' := \left\{\gamma_1,\gamma_2, \cdots , \gamma_n\right\}$.
 \item Let $L_0$ be the reduced lattice basis obtained by running $L(K;\Gamma')$ through LLL-reduction.
 \item Apply Lemma \ref{gram} along with the contrapositive of Lemma \ref{cond1} to $L_0$ to get a candidate for $N_\gamma$.
 \item For all $n < t \leq m$, let $L_t$ be the reduced lattice basis obtained from running $L\left(K;\Gamma' \cup \left\{\gamma_t\right\}\right)$ through LLL-reduction.
 \item Apply Lemma \ref{gram} with the contrapositive of Lemma \ref{cond2} to every $L_t$ to get $m-n$ additional candidates for $N_\gamma$.
 \item Set $N_\gamma$ to be the minimum of all of the candidates for $N_\gamma$ that were computed in steps 6 and 8.
\end{enumerate}
\caption{Returns an appropriate value for $N_\gamma$ that satisfies the conditions laid out in (\ref{ind1}) and (\ref{ind2}).}
\end{algorithm}
Note that by taking the minimum of all candidate values of $N_\gamma$, step 6 ensures that (\ref{ind1}) will be satisified. Moreover, step 6 ensures that (\ref{ind2}) will be satisfied when we consider $\gamma^* \in \Gamma'$ and step 8 ensures that (\ref{ind2}) will be satisfied when we consider $\gamma^* \not\in \Gamma'$.

\subsection{Finer Details}

\begin{defn}\label{LKGamma}
 Let $K \in \mathbb{Z}, S=\left\{\gamma_1,\cdots,\gamma_n\right\} \subset \mathbb{R}$ and define $\gamma_i'$ such that $K\gamma_i'=\lfloor K\gamma_i\rfloor$. Then $L(K;S) \subset \mathbb{Z}^{n+1}$ is the lattice generated by the following vectors

$$\left(\begin{array}{c}1\\ 0 \\ \vdots \\ 0 \\ K\gamma_1'\end{array}\right),
\left(\begin{array}{c}0\\ 1 \\ \vdots \\ 0 \\ K\gamma_2'\end{array}\right),
\cdots ,
\left(\begin{array}{c}0\\ \vdots \\ 0 \\ 1 \\ K\gamma_n'\end{array}\right).$$

\end{defn}

Our main result is centred around the following lattice since any vector in $L(K;S)$ will be of the form:
$$a_1\left(\begin{array}{c}1\\ 0 \\ \vdots \\ 0 \\ K\gamma_1'\end{array}\right)
+a_2\left(\begin{array}{c}0\\ 1 \\ \vdots \\ 0 \\ K\gamma_2'\end{array}\right)
+\cdots
+a_n\left(\begin{array}{c}0\\ \vdots \\ 0 \\ 1 \\ K\gamma_n'\end{array}\right)
= \left(\begin{array}{c}a_1\\ a_2 \\ \vdots \\ a_n \\ Kx\end{array}\right),$$ where $x = a_1\gamma_1'+a_2\gamma_2'+\cdots+a_n\gamma_n'$.

\begin{defn}\label{lin-dep}
 A set $S = \left\{\gamma_1,\gamma_2, \cdots, \gamma_n\right\} \subset \mathbb{R}$ is said to be $N_\gamma$-dependent if there exists $\alpha_1,\alpha_2,\cdots,\alpha_n \in \mathbb{R}$, not all zero, such that $\alpha_1\gamma_1+\cdots+\alpha_n\gamma_n= 0$ with $|\alpha_i| \leq N_\gamma$. $S$ is said to be $N_\gamma$-independent if it is not $N_\gamma$-dependent.
\end{defn}

\begin{defn}\label{weak-lin-dep}
 A set $S = \left\{\gamma_1,\gamma_2, \cdots, \gamma_n\right\} \subset \mathbb{R}$ is said to be weakly $N_\gamma$-dependent if there exists $\alpha_1,\alpha_2,\cdots,\alpha_n \in \mathbb{R}$, not all zero, such that $\alpha_1\gamma_1+\cdots+\alpha_n\gamma_n= 0$ with $|\alpha_i| \leq N_\gamma+1$ with at most one $i$ such that $|\alpha_i| = N_\gamma+1$. $S$ is said to be weakly $N_\gamma$-independent if it is not weakly $N_\gamma$-dependent.
\end{defn}

We wish to alert the reader to the potential confusion between $\{N_\gamma\}$-independent (from Definition \ref{defcon1})  and $N_\gamma$-independent (from Definition  \ref{lin-dep}).

\begin{Lem}\label{cond1}
If $S = \left\{\gamma_1,\gamma_2, \cdots, \gamma_n\right\} \subset \mathbb{R}$ is weakly $N_\gamma$-dependent, then there exists a nonzero vector $v\in L(K;S)$ such that $|v|^2 < (n^2+n)N_{\gamma}^2+(2n+2)N_{\gamma}+2$.

\begin{proof}
Consider the following vector $v \in L(K;S),$

$$v=\left(\begin{array}{c}\alpha_1\\ \alpha_2 \\ \vdots \\ \alpha_n \\ K(\alpha_1\gamma_1'+\cdots+\alpha_n\gamma_n')\end{array}\right).$$
The assumptions of the lemma show
\begin{equation*}
\begin{split}
|v|^2&=\alpha_i^2+\cdots+\alpha_n^2+K^2(\alpha_1\gamma_1'+\cdots+\alpha_n\gamma_n')^2 \\
&= \alpha_i^2+\cdots+\alpha_n^2+K^2\{(\alpha_1\gamma_1'+\cdots+\alpha_n\gamma_n')-(\alpha_1\gamma_1+\cdots+\alpha_n\gamma_n)\}^{2}, 
\end{split}
\end{equation*}
since $\alpha_1\gamma_1+\cdots+\alpha_n\gamma_n= 0$. Upon using the upper bounds on $|\alpha_{i}|$ and the fact that $|\gamma_{i} - \gamma_{i}'| < K^{-1}$, it follows that
\begin{equation*}
|v|^{2} < nN_{\gamma}^2+2N_{\gamma}+1+K^2\left(\frac{nN_{\gamma}+1}{K}\right)^2,
\end{equation*}
whence the lemma follows.
\end{proof}

\end{Lem}

Similarly, to account for the remaining zeroes, viz. $\gamma^*\notin\Gamma'$, we use

\begin{Lem}\label{cond2}
 If $S = \left\{\gamma_1,\gamma_2, \cdots, \gamma_n,\gamma_t\right\} \subset \mathbb{R}$ is $N_\gamma$-dependent where $\gamma_t \not\in\Gamma'$ is a zero, then there exists a nonzero vector $v\in L(K;S)$ such that $|v|^2 < (n^2+n)N_{\gamma}^2+2nN_{\gamma}+2$.

\begin{proof}
The proof follows that of Lemma \ref{cond1}; for each $i$, we have $|\alpha_{i}| \leq N_\gamma$.
\end{proof}
\end{Lem}

Note that in both lemmas above, the bounds are independent of our choice of $K$. The following lemma is true of all lattices.

\begin{Lem}\cite[Proposition 3.14]{small-vecs}\label{gram}
 Let $L \subset \mathbb{Z}^n$ be a lattice of dimension $m$. Let $\{b_i^*\}$ be the Gram--Schmidt orthogonalization of the basis of $L$. Then $|x|^2 \geq \min\left(|b_i^*|^2\right)$ for any nonzero $x \in L$. 
\end{Lem}

\begin{theorem}\label{cond3}
Let $L_0 = L(K;\Gamma')$ and let $L_t = L(K;\Gamma' \cup \{\gamma_t\})$ where $\gamma_t \in (\Gamma \cap [0,T]) \setminus \Gamma'$. We define $\left\{b_1,\cdots,b_n\right\}$ and $\left\{\beta_{t,1},\cdots,\beta_{t,n},\beta_{t,t}\right\}$ to be a basis for each lattice, respectively.

The elements of $\Gamma'$ are $\{N_{\gamma}\}$-independent in $\Gamma \cap [0, T]$ if 
$$\min\left(|b_i^*|^2\right) \geq (n^2+n)N_{\gamma}^2+(2n+2)N_{\gamma}+2$$
and
$$\min\left(|\beta_{t,i}^*|^2\right) \geq (n^2+n)N_{\gamma}^2+2nN_{\gamma}+2$$
for all $\gamma_t \in (\Gamma \cap [0,T]) \setminus \Gamma'$.

\begin{proof}
 We have two conditions to check. The first, (\ref{ind1}), is taken as a direct consequence of Lemma \ref{gram} and the contrapositive of Lemma \ref{cond1}. (\ref{ind2}) must be broken up into two separate parts. If $\gamma^* \in \Gamma'$, then we may apply the contrapositive of Lemma \ref{cond1} again. However, if $\gamma' \not\in \Gamma'$, then we must use the contrapositive of Lemma \ref{cond2}.
\end{proof}

\end{theorem}

Therefore, given a basis for our lattice, we may determine a lower bound for  $N_{\gamma}$. Note that we do not care which basis of the lattice we choose. At first glance, one may expect to take the basis given in Definition \ref{LKGamma}. However, if one attempts to perform the Gram--Schmidt on this basis, the vectors will be extremely {\it short}. It should only take a minute to convince the enthralled reader that even $|b_2^*|$ is relatively small. For this reason, we must find alternative bases for each lattice. Since there are many bases from which to choose we apply the LLL-reduction algorithm to find a {\it nearly orthogonal basis} for each lattice. By doing this, we will increase the length of the vectors obtained through the Gram--Schmidt orthogonalization process.

By definition of LLL-reduction, once the basis is reduced via LLL-reduction, it is guaranteed that $$|b_p^*|^2 \geq |b_{p-1}^*|^2\left(\delta-\frac14\right)$$ for all permissible $p$. These values are immensely important in our computation, since they give explicit information regarding the length of the shortest vector in the lattice. When computing the LLL-reduction, we tested several values of $\delta$ and $k$ to determine whether there was a significant difference in the choices. It turns out that the choice of $\delta$ is far less important than the choice of $k$. We refer the reader to \cite[Ch. 2]{computational-num-theory} for a closer look at LLL-reduction, including basically the same set up of vectors to determine the linear dependence of a set of numbers.

We chose $K = 10^k$ for some positive $k$: the $\gamma_i'$s are simply $\gamma_i$s accurate to $k$ decimal places (and then truncated). We used GP/Pari and Sage's functions that compute the zeroes. The programs were run independently, and we verified the zeroes using the intermediate value theorem on the Riemann $\xi$-function.

\subsection{An improved kernel}
In Theorem \ref{AndThm}, Anderson and Stark follow Ingham and use the F\'{e}jer kernel
\begin{equation*}
f(t) = \begin{cases} 1 - \frac{|t|}{T}, &\quad |t|\leq T, \\
0, &\quad |t| >T \end{cases}
\end{equation*}
to truncate the relevant sums. A permissible function for such an endeavour is one which has a non-negative Fourier transform and is supported on $[-T, T]$, and which is `close'  to unity in a neighbourhood about $t=0$. The last condition ensures that the contributions of the lower-lying zeroes are maximised. We use the function
\begin{equation*}
f_0(t) = \begin{cases} (1- \frac{|t|}{T})\cos \frac{\pi t}{T} + \pi^{-1}\sin\frac{\pi |t|}{T}, &\quad |t| \leq T, \\
0, &\quad |t|>T
\end{cases}
\end{equation*}
which is used by Odlyzko and te Riele in \cite{MertDis} --- for a discussion about the origin of this function see \cite[\S 4.1]{MertDis}.

\subsection{Sorting the zeroes}

Though the indices mentioned above may suggest that we must use the first $n$ zeroes as $\Gamma'$, this is not the case. Since we want to maximize equation (\ref{AndS3}), we shall sort the zeroes based on the ordering $\prec$ defined as follows
\begin{equation}\label{heavy-order}\gamma_i \prec \gamma_j \iff \frac{1}{|\rho_{i} \zeta'(\rho_{i})|} f_0\left(\gamma_{i}\right) > \frac{1}{|\rho_{j} \zeta'(\rho_{j})|} f_0\left(\gamma_{j}\right).\end{equation} We shall say that $\gamma_i$ is {\it heavier} than $\gamma_j$ if $\gamma_i \prec \gamma_j$.
For small values of $i$, sorting via $\prec$ does not seem to affect the order very much. However, as $T$ increases, the zeroes become scrambled.

\subsection{Results}

For our computation, we applied the main algorithm with $k = 9000$, $n = 500$ and $T \approx \gamma_{2001}-\epsilon$, where $\gamma_{2001}$ is the 2001st smallest zero. For the steps of the algorithm which required LLL-reduction, we used the standard $\delta = \frac34$ when performing LLL-reduction on $L_0$ and a weaker $\delta=\frac{3}{10}$ when reducing each $L_t$. Step 6 of the algorithm gave a candidate value for $N_\gamma$ of $794948$. When running the remaining 1500 zeroes through step 8, we find that the minimum candidate for $N_\gamma$ is $4976$. Thus, we arrive at the following

\begin{theorem}
Let $\Gamma'$ be the heaviest 500 zeroes with $T=\gamma_{2001}-\epsilon$. Then the elements of $\Gamma'$ are $\{N_\gamma\}$-independent, where all $N_\gamma$s are $4976$.
\end{theorem}

\subsection{Improvements}

Naturally, one should like to choose $\Gamma'$ to have as many entries as possible and $K$ to be as large as possible to improve on the bounds in Theorem \ref{Ours}. Unfortunately, the time taken to run each LLL-reduction is $O(n^{6}\log^3K) = O(n^6k^3)$, so either one of these choices may result in a quick computational explosion.

If we enlarge $n$ without changing $K$, the value of $N_\gamma$ is likely to decrease dramatically. Our experiments have shown that the value of $N_\gamma$ tends to drop to zero eventually as $n$ increases. On the flip-side, increasing the value of $K$ when the value of $N_\gamma$ is already large is relatively fruitless, as an increase of $K$ can only increase $N_\gamma$ and the $\frac{2N_\gamma}{N_\gamma+1}$ term is already close to $2$.

One might also suggest taking a larger value of $T$, re-sorting the zeroes and computing the corresponding $N_\gamma$. This seems to be the best possibility. By re-sorting the zeroes for each pair of $n$ and $T$, one has the {\it optimal} solution (provided the values of $N_\gamma$ remains large). However, if one wishes to roll the dice with different values of $n$ and $T$, one must do most of each computation from scratch. The first part of this recalculation can be cut down dramatically by storing specific intermediate results of the LLL-reduction and starting the reduction part way through. The second part of the recalculation, however, must be completely redone each time a new set of zeroes is selected. 

Say we fix $n$ and we wish to send $T \rightarrow \infty$, sorting the zeroes once again for each selection of $T$. In order for a high zero to have a large contribution, its derivative must be small. However, in \cite{numeric-deriv}, it is stated that {\it small} values of $|\zeta'(\frac12+i\gamma)|$, about 0.002, do not appear until $|\gamma| \approx 10^{22}$, meaning that their contribution to the sum will be minuscule. Thus, once $T$ is raised past a {\it reasonable} height it is unlikely that the first $n$ sorted zeroes will change.

Figure 1 shows the relationship between $n$ and $T$ with resorting of the zeroes. The chart assumes that $\frac{2N_\gamma}{N_\gamma+1} \approx 2$, which is a fair assumption if one believes that the zeroes are indeed linearly independent. Of special note, sorting the initial 9000 zeroes and taking the {\it best} 1000 in sorted order gives the first glimpse at improving Theorem \ref{Ours} by replacing 1.6383 with 2.

\begin{figure}\label{fig:nvst}\centering \includegraphics{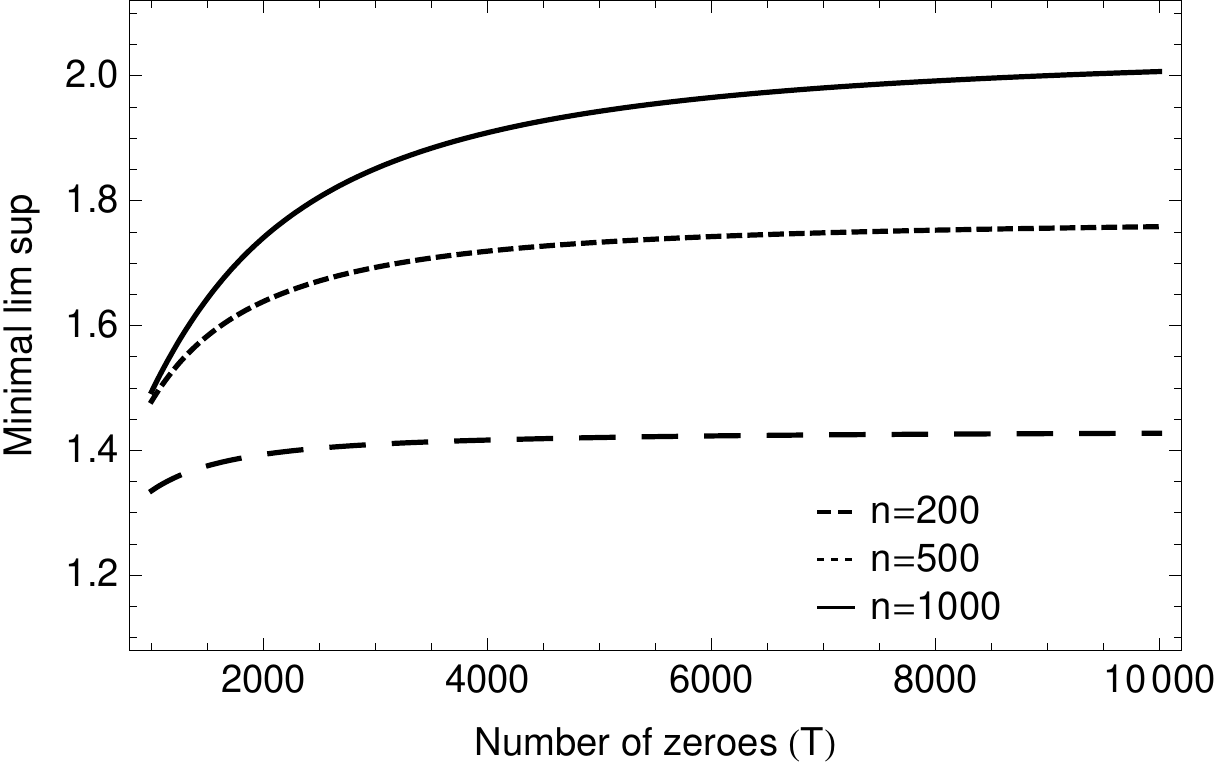}\caption{Smallest absolute value that the $\limsup$ and $\liminf$ can take if we use the first $n$  zeroes sorted by using the appropriate $T$, assuming $\frac{2N_{\gamma}}{N_{\gamma}+1} \approx 2$.}\end{figure}

To avoid the recomputation stated above, one may wish simply to  increase the value of $T$ without resorting the zeroes. Unfortunately all this will accomplish is making the kernel closer to 1, meaning we will eventually hit a ceiling. For illustrative purposes, Figure 2 shows the value Theorem \ref{Ours} could obtain if one were to raise the value of $T$ (assuming that the value of $N_\gamma$ stays large). It uses the first 300 zeroes in sorted order (sorted using $T = \gamma_{2001}-\epsilon$), but $T$ varies. We also include the values that the F\'{e}jer kernel would produce.

It is clear that the improved kernel approaches the maximal value quicker than the F\'{e}jer kernel. To obtain values within $0.01$ of the optimal value, the F\'{e}jer kernel needs to check a total of 30398 zeroes, while the improved kernel only needs to check 6224 zeroes. When a higher precision is needed, the gap widens immensely: to obtain values within $0.001$ of the optimal value, the F\'{e}jer kernel and the improved kernel need to check 399444 and 24043 zeroes, respectively.

\begin{figure}\label{fig:asymp}\centering \includegraphics{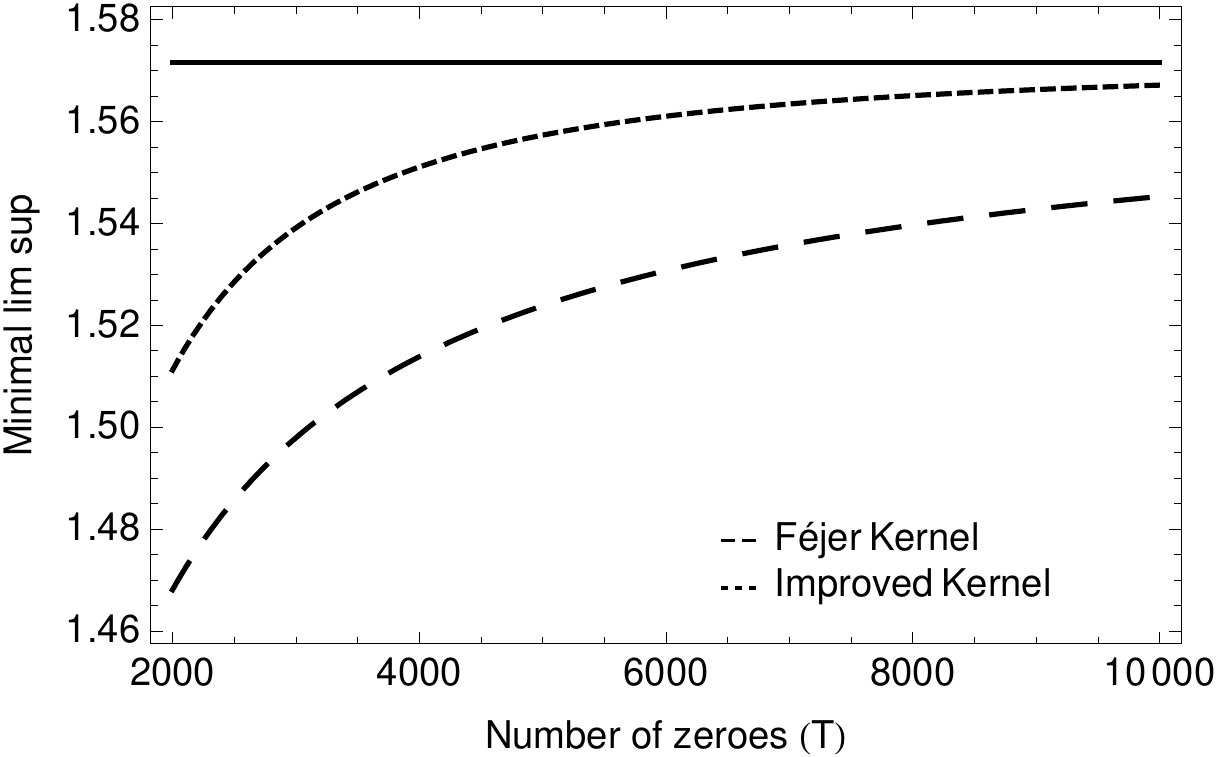}\caption{Smallest absolute value that the $\limsup$ and $\liminf$ can take if we use the first 300 zeroes sorted by using $T=\gamma_{2001}-\epsilon$ for the two different kernels. The horizontal line is the maximum attainable value for this set of zeroes.}\end{figure}

There is also a trade-off when determining which value to take for $\delta$ in the LLL-reduction. As we took more zeroes, the differences really started to shine through. Taking larger values of $\delta$ yielded a slower program, but one that gave a much better value for $N_{\gamma}$. On the other hand, taking a smaller value of $\delta$ sped up the program, but gave much smaller values for $N_{\gamma}$.

Finally, we have assumed that each of the $N_\gamma$s be the same. This is not necessary for our method to work. The bound in Lemma \ref{cond1} may be rewritten as $$|v|^2 <  2\max(N_\gamma) + 1 + \left(1+\sum N_{\gamma}\right)^2 + \sum N_{\gamma}^2.$$ A similar bound may be drawn up for Lemma \ref{cond2}. These bounds, however, are only useful when the resulting value of $N_{\gamma}$ from Theorem \ref{cond3} is relatively small.

\subsection{Other theorems}

In \cite{Grosswald}, two other important number-theoretic results are reproved on the condition that certain combinations of zeroes are linearly independent. Indeed, in \cite{Grosswald}, Theorem 1 gives us that $\pi(x) - \textrm{li}(x)$ changes sign infinitely often provided that the first 30 zeroes are $5$-independent; Theorem 2 shows that the functions associated with conjectures of P\'{o}lya and Tur\'{a}n, respectively
\begin{equation*}
L(x) = \sum_{1\leq n \leq x} \lambda(n), \quad T(x) = \sum_{1\leq n \leq x} \frac{\lambda(n)}{n},
\end{equation*}
where $\lambda(n)$ is the Liouville function, change sign infinitely often provided that the first 13 zeroes are $16$-independent. 

The data provided in Table \ref{n-indep} are more than enough to provide new proofs of these results.

\appendix
\section{$m$-Independence}

Table \ref{one-indep} gives the smallest sum (in absolute value) of the first $N$ zeroes of the zeta-function using coefficients $|c_n| \leq 1$. We have checked all permissible sums up to $N=41$. The first $20$ zeroes were checked in \cite[Table I]{Bateman}, labelled as {\it Type (A)}. They also provide a probabilistic value for the minimum value of such a sum, which we also include below.

To avoid a lengthy column of $\gamma$s to show the smallest linear combination, we encode the sums by an ordered pair of integers. If you write each integer in terms of its binary representation, a $1$ in the $i$th least significant bit implies that $\gamma_i$ is in the sum. The $i$th least significant bit being a $1$ in the first (resp. second) coordinate gives us a positive (resp. negative) coefficient. For example, $(5,24)$ represents the sum $\gamma_1+\gamma_3-\gamma_4-\gamma_5$.

\renewcommand{\arraystretch}{1.1}
\begin{table}[H]
\caption{Value of the Smallest Sums with Coefficients -1,0,1}
\centering
\begin{tabular}{@{}llll@{}}
\toprule \label{one-indep}
$n$ & Actual Value & Predicted Value & Linear Combination\\ \midrule
$20$ & $2.9799 \times 10^{-8}$  & $1.3976 \times 10^{-7}$  & $(533185,147768)$ \\
$21$ & $2.9799 \times 10^{-8}$  & $4.9104 \times 10^{-8}$  & $(533185,147768)$ \\
$22$ & $2.9799 \times 10^{-8}$  & $1.7238 \times 10^{-8}$  & $(533185,147768)$ \\
$23$ & $7.1672 \times 10^{-9}$  & $6.0341 \times 10^{-9}$  & $(3442980,4273746)$ \\
$24$ & $1.1632 \times 10^{-9}$  & $2.1088 \times 10^{-9}$  & $(2626459,12657764)$ \\
$25$ & $3.8873 \times 10^{-10}$ & $7.3493 \times 10^{-10}$ & $(17704982,10589760)$ \\
$26$ & $1.0788 \times 10^{-10}$ & $2.5605 \times 10^{-10}$ & $(42549638,3575905)$ \\
$27$ & $1.0788 \times 10^{-10}$ & $8.9049 \times 10^{-11}$ & $(42549638,3575905)$ \\
$28$ & $1.8340 \times 10^{-11}$ & $3.0897 \times 10^{-11}$ & $(96882844,171511617)$ \\
$29$ & $1.1519 \times 10^{-11}$ & $1.0713 \times 10^{-11}$ & $(93167683,405819176)$ \\
$30$ & $9.1777 \times 10^{-12}$ & $3.7102 \times 10^{-12}$ & $(948312448,41509390)$ \\
$31$ & $2.4115 \times 10^{-12}$ & $1.2836 \times 10^{-12}$ & $(1889619981,88484592)$ \\
$32$ & $4.6939 \times 10^{-14}$ & $4.4343 \times 10^{-13}$ & $(2299561158,1107850008)$ \\
$33$ & $4.6939 \times 10^{-14}$ & $1.5299 \times 10^{-13}$ & $(2299561158,1107850008)$ \\
$34$ & $4.6939 \times 10^{-14}$ & $5.2784 \times 10^{-14}$ & $(2299561158,1107850008)$ \\
$35$ & $1.8196 \times 10^{-17}$ & $1.8180 \times 10^{-14}$ & $(19757670928,14533859426)$ \\
$36$ & $1.8196 \times 10^{-17}$ & $6.2574 \times 10^{-15}$ & $(19757670928,14533859426)$ \\
$37$ & $1.8196 \times 10^{-17}$ & $2.1516 \times 10^{-15}$ & $(19757670928,14533859426)$ \\
$38$ & $1.8196 \times 10^{-17}$ & $7.3945 \times 10^{-16}$ & $(19757670928,14533859426)$ \\
$39$ & $1.8196 \times 10^{-17}$ & $2.5399 \times 10^{-16}$ & $(19757670928,14533859426)$ \\
$40$ & $1.8196 \times 10^{-17}$ & $8.7157 \times 10^{-17}$ & $(19757670928,14533859426)$ \\
$41$ & $1.8196 \times 10^{-17}$ & $2.9877 \times 10^{-17}$ & $(19757670928,14533859426)$ \\
\bottomrule
\end{tabular}
\end{table}

Table \ref{n-indep} gives us new lower bounds on the $m$-independence of the first $n$ zeroes of the zeta-function. By $m$-independence here, we mean that all non-trivial linear combinations of the first $n$ zeroes are nonzero assuming the coefficients are no more than $m$ in absolute value. This was computed using the same method as above, but keeping the zeroes in cardinal order. We include full results for the first 20 zeroes, and then only specific entries past there. Note that if $k$ is not included in the table, any bound for the first $n > k$ zeroes also gives a lower bound for the first $k$ zeroes.

\renewcommand{\arraystretch}{1.1}
\begin{table}[H]
\caption{The first $n$ zeroes of the zeta-function are $m$-independent}
\centering
\begin{tabular}{@{}llcll@{}}
\hline
\toprule
\multicolumn{1}{@{}c@{}}{$n$} & \multicolumn{1}{@{}c@{}}{$m$} && \multicolumn{1}{@{}c@{}}{$n$} & \multicolumn{1}{@{}c@{}}{$m$}\\
\cmidrule{1-2} \cmidrule{4-5}
 2 & $3.19683\times 10^{4499}$ &&  50 & $3.66786\times 10^{177}$\\
 3 & $7.01089\times 10^{2999}$ &&  75 & $6.96347\times 10^{116}$\\
 4 & $2.55333\times 10^{2249}$ && 100 & $1.83869\times 10^{86}$\\
 5 & $3.18071\times 10^{1799}$ && 125 & $4.96418\times 10^{67}$\\
 6 & $1.69018\times 10^{1499}$ && 150 & $1.90667\times 10^{55}$\\
 7 & $6.90883\times 10^{1284}$ && 175 & $1.35536\times 10^{46}$\\
 8 & $1.68884\times 10^{1124}$ && 200 & $1.13717\times 10^{39}$\\
 9 & $1.12832\times 10^{999}$  && 225 & $2.29079\times 10^{33}$\\
10 & $1.21351\times 10^{899}$  && 250 & $6.69056\times 10^{28}$\\
11 & $1.33521\times 10^{817}$  && 275 & $1.38130\times 10^{25}$\\
12 & $9.26711\times 10^{748}$  && 300 & $6.20938\times 10^{21}$\\
13 & $1.57289\times 10^{691}$  && 325 & $2.05342\times 10^{19}$\\
14 & $5.10452\times 10^{641}$  && 350 & $3.13279\times 10^{16}$\\
15 & $7.35106\times 10^{598}$  && 375 & $3.56683\times 10^{14}$\\
16 & $1.51957\times 10^{561}$  && 400 & $2.33172\times 10^{12}$\\
17 & $1.34818\times 10^{528}$  && 425 & $2.86453\times 10^{10}$\\
18 & $5.05309\times 10^{498}$  && 450 & $4.95180\times 10^{8}$\\
19 & $1.74671\times 10^{472}$  && 475 & $1.90299\times 10^{7}$\\
20 & $3.58761\times 10^{448}$  && 500 & $5.54632\times 10^{5}$\\
\bottomrule
 \end{tabular}
\label{n-indep}\end{table}

\section*{Acknowledgements}
We wish to thank Professor Nathan Ng for bringing our attention to this problem and Professor Soroosh Yazdani for pushing the idea of LLL-reduction into the picture. We also wish to thank Professors Brent, van de Lune, te Riele and Arias de Reyna for helping us to calculate the zeroes to high precision.

\bibliographystyle{plain}
\bibliography{themastercanada}

\end{document}